\documentclass[amsfonts]{article}
\usepackage{amssymb}

\newtheorem{theorem}{Theorem}[section]

\newtheorem{example}[theorem]{Example}
\newtheorem{proposition}[theorem]{Proposition}
\newtheorem{lemma}[theorem]{Lemma}

\newtheorem{remark}[theorem]{Remark}

\def\bR{\mathbb{R}}

\def\bC{\mathbb{C}}

\def\cB{\mathcal{B}}

\def\cE{\mathcal{E}}
\def\cF{\mathcal{F}}

\def\cH{\mathcal{H}}

\def\cP{\mathcal{P}}

\def\cS{\mathcal{S}}

\begin{document}

\title{The Stochastic Wave Equation
with Fractional Noise: a random field approach}

\author{Raluca M. Balan\footnote{Corresponding author. Department of Mathematics and Statistics, University of Ottawa,
585 King Edward Avenue, Ottawa, ON, K1N 6N5, Canada. E-mail address:
rbalan@uottawa.ca} \ \footnote{Research supported by a grant from
the Natural Sciences and Engineering Research Council of Canada.} \
and Ciprian A Tudor \footnote{Laboratoire Paul Painlev\'e,
Universit\'e de Lille 1,
 F-59655 Villeneuve d'Ascq, France. Email address: tudor@math.univ-lille1.fr. Associate member: SAMOS/MATISSE, Centre d'Economie de La
Sorbonne, Universit\'e de Panth\'eon-Sorbonne Paris 1, 90 rue de
Tolbiac, 75634 Paris Cedex 13, France. }}

\date{December 17, 2009}

\maketitle

\begin{abstract}
\noindent We consider the linear stochastic wave equation with
spatially homogenous Gaussian noise, which is fractional in time
with index $H>1/2$. We show that the necessary and sufficient
condition for the existence of the solution is a relaxation of the
condition obtained in \cite{dalang99}, when the noise is white in
time. Under this condition, we show that the solution is
$L^2(\Omega)$-continuous. Similar results are obtained for the heat
equation. Unlike the white noise case, the necessary and sufficient
condition for the existence of the solution in the case of the heat
equation is {\em different} (and more general) than the one obtained
for the wave equation.
\end{abstract}




\section{Introduction}

The random field approach to s.p.d.e.'s initiated in \cite{walsh86},
has become increasingly popular in the past few decades, as an
alternative to the semigroup approach developed in
\cite{daprato-zabczyk92}, or the analytic approach of
\cite{krylov99}.

Generally speaking, a random field solution of the (non-linear)
equation:
\begin{equation}
\label{spde}Lu(t,x)=\alpha(u(t,x))\dot W(t,x)+\beta(u(t,x)), \quad
t>0,x \in \bR^d
\end{equation} (with vanishing initial conditions) is a
collection $\{u(t,x), t \geq 0, x \in \bR^d\}$ of square integrable
random variables, which satisfy the following integral equation:
$$u(t,x)=\int_0^t
\int_{\bR^d}G(t-s,x-y)\alpha(u(s,y))W(ds,dy)+\int_0^t
\int_{\bR^d}G(t-s,x-y)\beta(u(s,y))dy ds,$$ provided that both
integrals above are well-defined (the first being a stochastic
integral). In this context, $L$ is a second-order partial
differential operator with constant coefficients, $G$ is the
fundamental solution of $Lu=0$, and $\dot W$ is a formal way of
denoting the random noise perturbing the equation.

When the equation is driven by a space-time white noise (i.e. a
Gaussian noise which has the covariance structure of a Brownian
motion in space-time), the random field solution exists only if the
spatial dimension is $d=1$. In this case, the stochastic integral
above is defined with respect to a martingale measure, and the
solution is well-understood for most operators $L$, in particular
for the heat and wave operators (see \cite{walsh86},
\cite{carmona-nualart88} or \cite{sanz-sole05}).

To obtain a random field solution in higher dimensions, one needs to
consider a different type of noise, which can be either Gaussian,
but with a spatially homogenous covariance structure given formally
by:
$$E[\dot W(t,x) \dot W(s,y)]=\delta(t-s)f(x-y),$$
or of Poisson type. Historically, the two approaches have been
initiated at about the same time (see \cite{mueller97},
\cite{dalang-frangos98}, \cite{millet-sanzsole99} for the wave
equation with Gaussian noise in dimension $d=2$, and
\cite{dalang-hou97}, \cite{saint-loubert98} for the Poisson case).

After the ingenious extension of the martingale measure stochastic
integral due to \cite{dalang99}, it became clear that the random
field approach can be pursued for the study of general s.p.d.e.'s
with spatially homogenous Gaussian noise. Since this extension
allows for integrands which are non-negative measures (in space),
the theory developed in \cite{dalang99} covers instantly the case of
the (non-linear) wave equation in dimensions $d \in \{1,2,3\}$, and
the case of the heat equation in any dimensions $d$. In the
non-linear case, the existence of the solution is obtained by a
Picard's iteration scheme, under the usual Lipschitz assumptions on
$\alpha, \beta$, and the following condition, linking the operator
$L$ and the spatial covariance function $f$:
\begin{equation}
\label{Dalang-cond}\int_{\bR^d} \int_0^t |\cF G(u,\cdot)(\xi)|^2 du
\mu(d\xi)<\infty. \end{equation} (Here $\mu$ is a non-negative
tempered measure, whose Fourier transform in $f$.)

Moreover, (\ref{Dalang-cond}) is the necessary any sufficient
condition for the stochastic integral $\int_0^t
\int_{\bR^d}G(t-s,x-y)W(ds,dy)$ to be well-defined, and hence the
necessary any sufficient condition for the existence of the solution
in the linear case, when $\alpha \equiv 1$ and $\beta \equiv 0$.
Since for both heat and wave operators,
\begin{equation}
\label{estimates} c_t^{(1)}\frac{1}{1+|\xi|^2} \leq \int_0^t |\cF
G(u,\cdot)(\xi)|^2 du \leq c_t^{(2)}\frac{1}{1+|\xi|^2}, \quad
\mbox{for all} \ \xi \in \bR^d,
\end{equation} for some constants
$c_t^{(1)},c_t^{(2)}>0$, condition (\ref{Dalang-cond}) is equivalent
to:
$$\int_{\bR^d} \frac{1}{1+|\xi|^2}\mu(d\xi)<\infty.$$

Subsequently, using the Malliavin calculus techniques, it was shown
that the random variable $u(t,x)$ has an absolutely continuous law
with respect to the Lebesgue measure on $\bR$, and this density is
infinitely differentiable. These results are valid for the heat
equation in any dimension $d$, and for the wave equation in
dimension $d \in \{1,2,3\}$ (see \cite{QS-SanzSole04a},
\cite{QS-SanzSole04b}, \cite{sanz-sole05}), under the additional
assumption (which was removed in \cite{Nualart-QS07}):
\begin{equation}
\label{Dalang-cond-alpha} \int_{\bR^d} \left(\frac{1}{1+|\xi|^2}
\right)^{\alpha}\mu(d\xi)<\infty, \quad \mbox{for some} \ \alpha \in
(0,1).
\end{equation}
Under (\ref{Dalang-cond-alpha}), one also obtains the H\"{o}lder
continuity of the solution for the heat equation in any dimension
$d$ and the wave equation in dimensions $d \in \{1,2,3\}$. This is
done using Kolmogorov's criterion and some estimates for the $p$-th
moments of the increments of the solution (see
\cite{sanzsole-sarra00}, \cite{sanzsole-sarra02},
\cite{dalang-sanzsole08}).

The case of the wave equation in dimension $d \geq 4$ was solved in
the recent article \cite{conus-dalang09}, using an extension of the
integral developed in \cite{dalang99}. The existence of a
random-field solution is obtained under condition
(\ref{Dalang-cond}). In the the affine case (i.e. $\alpha(u)=au+b,
a,b \in \bR$ and $\beta \equiv 0$), and under the additional
assumption (\ref{Dalang-cond-alpha}), the solution is shown to be
H\"older continuous.

In parallel with these developments, a new process began to be used
intensively in stochastic analysis: the {\em fractional Brownian
motion} (fBm) with index $H \in (0,1)$, a zero-mean Gaussian process
$(B_t)_{t \geq 0}$ with covariance:
$$R_{H}(t,s)=\frac{1}{2}(t^{2H}+s^{2H}-|t-s|^{2H}).$$

The case $H=1/2$ corresponds to the classical Brownian motion, while
the cases $H>1/2$ and $H<1/2$ have many contrasting properties. We
refer the reader to the survey article \cite{nualart03} and the
monographs \cite{BHOZ08} and \cite{mishura08} for more details. Most
importantly, in the case $H>1/2$,
\begin{equation}
\label{formula-RH} R_{H}(t,s)=\alpha_H\int_0^t
\int_0^s|u-v|^{2H-2}dudv,
\end{equation}
where $\alpha_H=H(2H-1)$. This shows that $(B_t)_{t \geq 0}$ has a
homogenous covariance structure, similar to the spatial structure of
the noise $\dot W$ considered above.

Returning to our discussion about s.p.d.e.'s with a Gaussian noise,
it seems natural to consider equation (\ref{spde}), when the
covariance of the noise $\dot W$ is given formally by:
\begin{equation}
\label{noise} E[\dot W(t,x) \dot W(s,y)]=\alpha_H|t-s|^{2H-2}f(x-y).
\end{equation}

However, this simple modification changes the problem drastically,
since unless $H=1/2$, the fBm is {\em not} a semimartingale, and
therefore the previous method, based on martingale measure
stochastic integrals, cannot be applied.

Several methods have been proposed for developing a stochastic
calculus with respect to fBm: (i) the Malliavin calculus (see
\cite{decreusefond-ustunel99}, \cite{AMN01}, \cite{alos-nualart03},
\cite{nualart06}), which exploits the fact that the fBm is Gaussian;
(ii) the method of generalized Lebesgue-Stieltjes integration (see
\cite{zahle98}), which uses the H\"older continuity of the fBm
trajectories; (iii) the rough path analysis (see \cite{lyons98},
\cite{lyons-qian02}), which uses the fact that the paths of the fBm
have bounded $p$-variation, for $p>1/H$; (iv) the stochastic
calculus via regularization based also in general on the  properties
of the paths of the fBm (see \cite{GRV}).

These methods have been applied to s.p.d.e.'s (see
\cite{maslowski-nualart03}, \cite{nualart-vuillermont06},
\cite{sanzsole-vuillermont07}, \cite{GLT06}), \cite{QS-tindel07}),
but not using the random field approach
A notable exception is the heat equation. The linear equation with
noise (\ref{noise}) and $H>1/2$ was examined in
\cite{balan-tudor08}, for particular functions $f$ (e.g.
$f(x)=|x|^{-(d-\alpha)}$ with $\alpha \in (0,d)$). We also mention
the works \cite{EV} and \cite{TTV} for the case of the space
variable belonging to the unit circle. The quasi-linear equation
(i.e. $\alpha \equiv 0$) was treated in \cite{oksendal-zhang01}, and
the equation with multiplicative noise (i.e. $\alpha(u)=u, \beta
\equiv 0)$ was studied in \cite{hu01}; in these two references, the
covariance structure of the noise is a particular case of
(\ref{noise}): for $H,H_i>1/2$
 $$E[\dot W(t,x) \dot
W(s,y)]=\alpha_H|t-s|^{2H-2}\prod_{i=1}^{d}(\alpha_{H_i}|x_i-y_i|^{2H_i-2}).$$
(This type of noise is called {\em fractional Brownian field}.)
The heat equation with multiplicative noise (\ref{noise}) was
studied in \cite{balan-tudor09} (for particular functions $f$ and
$H>1/2$) and \cite{hu-nualart09} (in the case $H \in (0,1)$ and
$f=\delta_0$). In the  case when the spatial dimension is $d=1$, the
non-linear equation has been treated in \cite{QS-tindel07} using a
two-parameter Young integral based on the H\"older continuity of
fBm.

To the best of our knowledge, there is no study of the wave equation
driven by a noise $\dot W$, whose covariance is given by
(\ref{noise}). The goal of the present article is to start filling
this gap, by identifying the necessary and sufficient conditions for
the existence of a random field solution of the linear wave equation
with noise (\ref{noise}) and $H>1/2$. We also treat the heat
equation.

When $H>1/2$, it turns out that under relatively mild assumptions on
the fundamental solution $G$ of the operator $L$, the necessary and
sufficient condition for the existence of the random-field solution
of the linear equation $Lu=\dot W$ is:
\begin{equation}
\label{general-cond} \int_{\bR^d}\int_0^t \int_0^t \cF
G(u,\cdot)(\xi)
 \overline{\cF G(v,\cdot)(\xi)} |u-v|^{2H-2} dudv  \mu(d\xi)<\infty,
\end{equation}
which is more general than (\ref{Dalang-cond}). Note that the
integrand of the $\mu(d\xi)$ integral in (\ref{general-cond}) is the
$\cH(0,t)$-norm of the function $u \mapsto \cF G(u,\cdot)(\xi)$.
Quite surprisingly, and in contrast with (\ref{estimates}), the
estimates that we obtain for this norm are {\em different} in the
case of the wave and heat operators: in the case of the wave
equation, (\ref{general-cond}) is equivalent to
\begin{equation}
\label{wave-cond} \int_{\bR^d}
\left(\frac{1}{1+|\xi|^2}\right)^{H+1/2}\mu(d\xi)<\infty,
\end{equation}
whereas in the case of the heat equation, (\ref{general-cond}) is
equivalent to:
\begin{equation}
\label{heat-cond} \int_{\bR^d}
\left(\frac{1}{1+|\xi|^2}\right)^{2H}\mu(d\xi)<\infty,
\end{equation}

The amazing fact is that for the wave operator, these estimates can
be deduced using only the estimates of the $L^2(0,t)$-norm (given by
(\ref{estimates})), the trick being to pass to the spectral
representation of the $\cH(0,t)$-norm of $u \mapsto \cF
G(u,\cdot)(\xi)$. In the case of the heat operator, there is no need
for this machinery, since $u \mapsto \cF G(u,\cdot)(\xi)$ is a
non-negative function, and its $\cH(0,t)$-norm can be bounded
directly by the $L^{1/H}(0,t)$-norm, which is easily computable.

This article is organized as follows. Section 2 contains some
preliminaries, and a basic result which ensures that under
(\ref{general-cond}), the stochastic integral of the fundamental
solution $G$ of the wave operator is well defined. In Section 3, we
show that the solution of the wave equation exists if and only if
(\ref{wave-cond}) holds (Theorem \ref{wave-th}). Moreover, the
solution is $L^2(\Omega)$-continuous. Similar results are obtained
in Section 4 for the heat equation, using (\ref{heat-cond}).
Appendix A contains some useful identities, which are needed in the
sequel. Appendix B gives the spectral representation of the
$\cH(0,t)$-norm of the function $\sin$.

\section{The Basics}

We denote by $C_0^{\infty}(\bR^{d+1})$ the space of infinitely
differentiable functions on $\bR^{d+1}$ with compact support, and
$\cS(\bR^d)$ the Schwartz space of rapidly decreasing $C^{\infty}$
functions in $\bR^d$. For $\varphi \in L^1(\bR^d)$, we let $\cF
\varphi$ be the Fourier transform of $\varphi$:
$$\cF \varphi (\xi)=\int_{\bR^d} e^{-i \xi \cdot x}\varphi (x)dx.$$

We begin by introducing the framework of \cite{dalang99}. 
Let $\mu$ be a non-negative tempered measure on $\bR^d$, i.e. a
non-negative measure which satisfies: $$\int_{\bR^d}
\left(\frac{1}{1+|\xi|^2} \right)^l \mu(d\xi)<\infty, \quad
 \mbox{for some} \ l >0.$$

Since the integrand is non-increasing in $l$, we may assume that $l
\geq 1$ is an integer. Note that $1+|\xi|^2$ behaves as a constant
around $0$, and as $|\xi|^2$ at $\infty$, and hence
(\ref{mu-tempered}) is equivalent to:
\begin{equation}
\label{mu-tempered} \int_{|\xi| \leq 1}\mu(d\xi)<\infty, \quad
\mbox{and} \quad \int_{|\xi| \geq 1}\frac{1}{|\xi|^{2l}}<\infty,
\quad
 \mbox{for some integer} \ l \geq 1.
  \end{equation}

Let $f: \bR^d \to \bR_{+}$ be the Fourier transform of $\mu$ in
$\cS'(\bR^d)$, i.e.
$$\int_{\bR^d}f(x)\varphi(x)dx=\int_{\bR^d}\cF
\varphi(\xi)\mu(d\xi), \quad \forall \varphi \in \cS(\bR^d).$$

Simple properties of the Fourier transform show that for any
$\varphi, \psi \in \cS(\bR^d)$,
$$\int_{\bR^d} \int_{\bR^d} \varphi(x)f(x-y)\psi(y)dx dy=
\int_{\bR^d}\cF \varphi(\xi) \overline{\cF \psi(\xi)}\mu(d\xi).$$

An approximation argument shows that the previous equality also
holds for indicator functions $\varphi=1_{A},\psi=1_{B}$, with $A,B
\in \cB_{b}(\bR^d)$, where $\cB_b(\bR^d)$ is the class of bounded
Borel sets of $\bR^d$:
\begin{equation}
\label{Fourier-indicator} \int_A \int_B f(x-y)dx dy=\int_{\bR^d}\cF
1_{A}(\xi) \overline{\cF 1_{B}(\xi)} \mu(d\xi).
\end{equation}

As in \cite{balan-tudor08}, \cite{balan-tudor09}, on a complete
probability space $(\Omega,\cF,P)$, we consider a zero-mean Gaussian
process $W=\{W_t(A); t \geq 0, A \in \cB_{b}(\bR^d)\}$ with
covariance:
$$E(W_t(A)W_s(B))=R_{H}(t,s) \int_{A} \int_{B} f(x-y)dx dy=:
\langle 1_{[0,t] \times A}, 1_{[0,s] \times B} \rangle_{\cH \cP}.$$

Let $\cE$ be the set of linear combinations of elementary functions
$1_{[0,t] \times A}$, $t \geq 0, A \in \cB_b(\bR^d)$, and $\cH \cP$
be the Hilbert space defined as the closure of $\cE$ with respect to
the inner product $\langle \cdot , \cdot \rangle_{\cH \cP}$.
(Alternatively, $\cH \cP$ can be defined as the completion of
 $C_0^{\infty}(\bR^{d+1})$, with respect to the inner product
$\langle \cdot, \cdot \rangle_{\cH \cP}$.)

The map $1_{[0,t] \times A} \mapsto W_t(A)$ is an isometry between
$\cE$ and the Gaussian space $H^{W}$ of $W$, which can be extended
to $\cH \cP$. We denote this extension by: $$\varphi \mapsto
W(\varphi)=\int_0^{\infty}\int_{\bR^d} \varphi(t,x)W(dt,dx).$$

In the present work, we assume that $H>1/2$. Hence,
(\ref{formula-RH}) holds. From (\ref{Fourier-indicator}) and
(\ref{formula-RH}), it follows that for any $\varphi,\psi \in \cE$,
\begin{eqnarray*}
\langle \varphi, \psi \rangle_{\cH \cP}&=& \alpha_H
\int_{0}^{\infty} \int_0^{\infty}\int_{\bR^d}\int_{\bR^d}
\varphi(u,x)\psi(v,y)f(x-y)|u-v|^{2H-2} dx
dy du dv  \\
&=& \alpha_H  \int_{0}^{\infty} \int_0^{\infty} \int_{\bR^d} \cF
\varphi(u,\cdot)(\xi)\overline{\cF\psi(v,\cdot)(\xi)} |u-v|^{2H-2}
\mu(d\xi) du dv.
\end{eqnarray*}

Moreover, we can interchange the order of the integrals $dudv$ and
$\mu(d\xi)$, since for indicator functions $\varphi$ and $\psi$, the
integrand is a product of a function of $(u,v)$ and a function of
$\xi$. Hence, for $\varphi,\psi \in \cE$, we have:
\begin{equation}
\label{norm-HP-2} \langle \varphi, \psi \rangle_{\cH \cP}= \alpha_H
\int_{\bR^d} \int_{0}^{\infty} \int_0^{\infty} \cF
\varphi(u,\cdot)(\xi)\overline{\cF\psi(v,\cdot)(\xi)} |u-v|^{2H-2}
du dv \mu(d\xi).
\end{equation}


The space $\cH \cP$ may contain distributions, but contains the
space $|\cH \cP|$ of measurable functions $\varphi: \bR_{+} \times
\bR^d \to \bR$ such that
$$\|\varphi \|_{|\cH \cP|}^2:=\alpha_H
\int_{0}^{\infty} \int_0^{\infty}\int_{\bR^d}\int_{\bR^d}
|\varphi(u,x)||\varphi(v,y)|f(x-y)|u-v|^{2H-2} dx dy du dv<\infty.$$


We recall now several facts related to the fBm (see e.g.
\cite{nualart03}).

Let $B=(B_t)_{t \geq 0}$ be a fBm of index $H>1/2$. For a fixed
$T>0$, let $\cH(0,T)$ be the Hilbert space defined as the closure of
$\cE(0,T)$ (the set of step functions on $[0,T]$), with respect to
the inner product:
$$\langle 1_{[0,t]}, 1_{[0,s]} \rangle_{\cH(0,T)}=R_{H}(t,s).$$

One can prove that
$$R_{H}(t,s)=\int_{0}^{t \wedge s} K_{H}(t,r)K_H(s,r)dr,$$
where $K_H(t,r)= c_H^*\int_r^t (u-r)^{H-3/2} u^{H-1/2}du$ and
$c_H^*=\left(\frac{\alpha_H}{\beta(H-1/2,2-2H)}\right)^{1/2}$. (Here
$\beta$ denotes the Beta function.) Therefore, the map $K_H^*$
defined by:
$$(K_H^* 1_{[0,t]}) (s)=K_{H}(t,s)1_{[0,t]}(s)$$
is an isometry between $\cE(0,T)$ and $L^2(0,T)$. This isometry can
be extended to $\cH(0,T)$, and is denoted by $\phi \mapsto
B(\phi)=\int_0^T \phi(s)dB_s$.

The transfer operator $K_H^*$ can be expressed in terms of
fractional integrals, as follows: for any $\phi \in \cE(0,T)$,
$$(K_H^* \phi)(s)=c_{H}^{*}\Gamma(H-1/2)
s^{1/2-H}I_{T-}^{H-1/2}(u^{H-1/2}\phi(u))(s),$$ where
$$I_{T-}^{\alpha}f (s)=
\frac{1}{\Gamma(\alpha)}\int_{s}^{T}(u-s)^{\alpha-1}f(u)du$$ denotes
the fractional integral of $f \in L^1(0,T)$, of order $\alpha \in
(0,1)$.

$K_H^*$ can be extended to complex-valued functions, as follows. Let
$\cE_{\bC}(0,T)$ be the set of all complex linear combinations of
functions $1_{[0,t]}, t \in [0,T]$, and $\cH_{\bC}(0,T)$ be the
closure of $\cE_{\bC}(0,T)$ with respect to the inner product:
$$\langle \varphi,\psi \rangle_{\cH_{\bC}(0,T)}=
\alpha_H \int_0^T \int_0^T \varphi(u)
\overline{\psi(v)}|u-v|^{2H-2}du dv.$$

The operator $K_H^*$ is an isometry which maps $\cH_{\bC}(0,T)$ onto
$L_{\bC}^2(0,T)$ (the space of functions $\varphi:[0,T] \to \bC$,
with $\int_0^T |\varphi(t)|^2 dt<\infty$): for any $\phi \in
\cH_{\bC}(0,T)$,
\begin{equation}
\label{isometry} \alpha_H \int_0^T \int_0^T \phi(u)
\overline{\phi(v)}|u-v|^{2H-2}du dv=d_{H} \int_0^T
|I_{T-}^{H-1/2}(u^{H-1/2}\phi(u))(s)|^2 \lambda_H(ds),
\end{equation} where $d_H=(c_{H}^{*})^2\Gamma(H-1/2)^2$ and
$\lambda_H(ds)=s^{1-2H}ds$.

Let $\cE_T$ be the class of elementary functions on $[0,T] \times
\bR^d$.  Note that for any $\varphi \in \cE_T$, the function $t
\mapsto \cF \varphi (t,\cdot)(\xi)$ belongs to $\cH_{\bC}(0,T)$, for
all $\xi \in \bR^d$. Using (\ref{norm-HP-2}) and (\ref{isometry}),
we obtain that for any $\varphi \in \cE_T$,
\begin{equation}
\label{new-def-norm-HP} \|\varphi \|_{\cH \cP}^2=d_H \int_{\bR^d}
\int_0^T |I_{T-}^{H-1/2}(u^{H-1/2} \cF \varphi (u,\cdot)(\xi))(s)|^2
\lambda_H(ds) \mu(d\xi)=:\|\varphi\|_{0}^2.
\end{equation}

We are now ready to state our result. Note that, although the
conclusion of this result resembles that of Theorem 3 of
\cite{dalang99} (for deterministic integrands), the hypothesis are
different, since the proof uses techniques specific to the fBm.

\begin{theorem}
\label{theorem-about-HP} Let $[0,T] \ni t \mapsto \varphi(t,\cdot)
\in \cS'(\bR^d)$ be a deterministic function such that $\cF
\varphi(t,\cdot)$ is a function for all $t \in [0,T]$. Suppose that:\\
 (i) the function $t \mapsto \cF \varphi(t,
\cdot)(\xi)$ belongs to
$\cH_{\bC}(0,T)$ for all $\xi \in \bR^d$; \\
(ii) the function $(t,\xi) \mapsto \cF \varphi(t, \cdot)(\xi)$ is
measurable on $(0,T) \times \bR^d$; \\
(iii) $\int_s^T u^{H-1/2} (u-s)^{H-3/2} |\cF
\varphi(u,\cdot)(\xi)|du<\infty$ for all $(s,\xi) \in (0,T) \times
\bR^d$ (or $\cF \varphi (s,\cdot)(\xi) \geq 0$ for all $(s,\xi) \in
(0,T) \times \bR^d$).

If
\begin{equation}
\label{HP-norm-finite}I_T:=\alpha_H \int_{\bR^d} \int_0^T \int_0^T
\cF \varphi(u,\cdot)(\xi) \overline{\cF
\varphi(v,\cdot)(\xi)}|u-v|^{2H-2}du dv \mu(d\xi)<\infty,
\end{equation}
then $\varphi \in \cH \cP$ and $\|\varphi \|_{\cH \cP}^2 =I_T$. (By
convention, we set $\varphi(t,\cdot)=0$ for $t>T$.)
\end{theorem}

\begin{remark}
{\rm Conditions (i)-(iii) are satisfied by the fundamental solution
$G$ of the wave (or heat) equation. In this case, $|\cF
G(t,\cdot)(\xi)| \leq 1$ for all $t \geq 0$, and hence the map $t
\mapsto \cF G(t,\cdot)(\xi)$ belongs to $L_{\bC}^2(0,T)$, which is
included in $\cH_{\bC}(0,T)$.}
\end{remark}

\noindent {\bf Proof:} The argument is a modified version of the
proof of Theorem 3.8 of \cite{balan-tudor08}. For any $\xi \in
\bR^d$ fixed, we apply (\ref{isometry}) to the function
$\phi_{\xi}(t)=\cF \varphi(t,\cdot)(\xi)$. We get:
\begin{eqnarray}
\label{equality-xi} \lefteqn{ \alpha_H \int_0^T \int_0^T \cF
\varphi(u,\cdot)(\xi) \overline{\cF
\varphi(v,\cdot)(\xi)}|u-v|^{2H-2}du dv=} \\
\nonumber  & & d_H \int_0^T|I_{T-}^{H-1/2}(u^{H-1/2}\cF
\varphi(u,\cdot)(\xi))(s)|^2 \lambda_H(ds).
\end{eqnarray}

It will be shown later that:
\begin{equation}
\label{a-measurable} (s,\xi) \mapsto
a(s,\xi):=I_{T-}^{H-1/2}(u^{H-1/2}\cF\varphi(u,\cdot)(\xi))(s)\
\mbox{is measurable on} \ (0,T) \times \bR^d.
\end{equation}

\noindent Hence, we can integrate with respect to $\mu(d\xi)$ in
(\ref{equality-xi}). Using (\ref{HP-norm-finite}), we obtain:
\begin{equation}
\label{a-in-L2} I_T=d_{H}\int_{\bR^d} \int_0^T
|I_{T-}^{H-1/2}(u^{H-1/2}\cF \varphi(u,\cdot)(\xi))(s)|^2
\lambda_H(ds) \mu(d\xi)=:\|\varphi \|_{0}^{2}<\infty.
\end{equation}

By the definition of $\cH \cP$ and (\ref{new-def-norm-HP}), it
suffices to show that for any $\varepsilon>0$, there exists a
function $l=l_{\varepsilon} \in \cE_T$ such that:
\begin{equation}
\label{E-is-dense} \|\varphi -l \|_{0} <\varepsilon.
\end{equation}


Let $\varepsilon>0$ be arbitrary. By (\ref{a-measurable}) and
(\ref{a-in-L2}), it follows that $a \in L^2((0,T) \times \bR^d,
\lambda_H(ds) \times \mu(d\xi) )$. Hence, there exists a simple
function $h(s,\xi)$ such that
\begin{equation}
\label{approx-a-h} \int_{\bR^d}
\int_{0}^{T}|a(s,\xi)-h(s,\xi)|^{2}\lambda_{H}(ds)\mu(d\xi)
<\varepsilon.
\end{equation}

Without loss of generality, we assume that
$h(s,\xi)=1_{(c,d]}(s)1_{A}(\xi)$, with $c,d \in [0,T],c < d$ and $A
\in \cB_{b}(\bR^d)$. By relation (8.1) of \cite{pipiras-taqqu01}, we
approximate the function $1_{(c,d]}(s)$ in
$L^2((0,T),\lambda_{H}(ds))$ by $I_{T-}^{H-1/2}(u^{H-1/2}l_0(u))(s)$
with $l_0 \in \cE(0,T)$, i.e.
\begin{equation}
\label{approx-l0} \int_0^T
|1_{(c,d]}(s)-I_{T-}^{H-1/2}(u^{H-1/2}l_0(u))(s)|^2\lambda_{H}(ds)
<\varepsilon.
\end{equation}

By Lemma 3.7 of \cite{balan-tudor08}, we approximate the function
$1_{A}(\xi)$ in $L^2(\bR^d,\mu(d\xi))$ by $\cF l_1(\xi)$ with $l_1
\in \cE(\bR^d)$, i.e.
\begin{equation}
\label{approx-l1} \int_{\bR^d}|1_{A}(\xi)-\cF
l_1(\xi)|^2\mu(d\xi)<\varepsilon.
\end{equation}

We define $l(u,x)=l_0(u)l_1(x)$. Clearly $l \in \cE_T$ and $\cF
l(u,\cdot)(\xi)=l_0(u)\cF l_{1}(\xi)$. Let
$$b(s,\xi):=I_{T-}^{H-1/2}(u^{H-1/2} \cF l (u,\cdot)(\xi))(s)=
I_{T-}^{H-1/2}(u^{H-1/2} l_0(u))(s) \cdot \cF l_1(\xi).$$

\noindent Using (\ref{approx-l0}) and (\ref{approx-l1}), we obtain
that:
\begin{eqnarray}
\nonumber \lefteqn{\int_{\bR^d}\int_0^T |h(s,\xi)-b(s,\xi)|^2
\lambda_{H}(ds)\mu(d\xi) }
\\
\nonumber & \leq & 2 \left\{ \int_{\bR^d}\int_0^T
|1_{(c,d]}(s)-I_{T-}^{H-1/2}(u^{H-1/2}
l_0(u))(s)|^2 1_{A}(\xi) \lambda_{H}(ds)\mu(d\xi) + \right. \\
\nonumber & & \left. \int_{\bR^d}\int_0^T |I_{T-}^{H-1/2}(u^{H-1/2}
l_0(u))(s)|^2
|1_{A}(\xi)-\cF l_1(\xi)|^2 \lambda_{H}(ds)\mu(d\xi) \right\}\\
\label{approx-b-h}  & \leq & 2
\{\varepsilon\mu(A)+\varepsilon\|l_0\|_{\cH
(0,T)}^2/d_H\}:=C_1\varepsilon.
\end{eqnarray}

From (\ref{approx-a-h}) and (\ref{approx-b-h}), it follows that
$$\|\varphi-l
\|_{0}^{2}=d_H\int_{\bR^d}
\int_{0}^{T}|a(s,\xi)-b(s,\xi)|^{2}\lambda_{H}(ds)\mu(d\xi) <2d_H
(\varepsilon+C_1\varepsilon):=C_2 \varepsilon.$$ This concludes the
proof of (\ref{E-is-dense}).

We now return to the proof of (\ref{a-measurable}), which uses
assumptions (ii) and (iii). If $\cF (u,\cdot)(\xi) \geq 0$, then
$(u,s,\xi) \mapsto \phi(u,s,\xi) =1_{\{s \leq
u\}}u^{H-1/2}(u-s)^{H-3/2}\cF \varphi(u,\cdot)(\xi)$ is measurable
and non-negative, and $a(s,\xi)=\int_0^T \phi(u,s,\xi)du$ is
measurable, by Fubini's theorem.

Suppose next that $\int_s^T u^{H-1/2}(u-s)^{H-3/2}|\cF \varphi
(u,\cdot)(\xi)|du<\infty$. If $l(s,\xi)=1_{(c,d]}(s)1_{A}(\xi)$ is
an elementary function with $c,d \in [0,T],A \in \cB_b(\bR^d)$, then
$$a_{l}(s,\xi)=I_{T-}^{H-1/2}(u^{H-1/2}l(u,\xi))(s)=1_{A}(\xi)\int_s^T
u^{H-1/2}1_{(c,d]}(u)(u-s)^{H-3/2}du$$ is clearly measurable. 
In general, since $(u,\xi) \mapsto \cF \varphi(u,\cdot)(\xi)$ is
measurable, there exists a sequence $(l_n)_n$ of simple functions
such that $l_n(u,\xi) \to \cF \varphi(u,\cdot)(\xi)$ for all
$(u,\xi)$ and $|l_n(u,\xi)| \leq |\cF \varphi(u,\cdot)(\xi)|$ for
all $(u,\xi),n$ (see e.g. Theorem 13.5 of \cite{billingsley95}). By
the dominated convergence theorem, for every $(s,\xi)$
$$|a_{l_n}(s,\xi)-a(s,\xi)| \leq \int_s^T u^{H-1/2} (u-s)^{H-3/2}|l_n(s,\xi)-
\cF \varphi (u,\cdot)(\xi)|du \to 0.$$ Since $a_{l_n}(s,\xi)$ is
measurable for every $n$, it follows that $a(s,\xi)$ is measurable.
$\Box$

\section{The wave equation}

We consider the linear wave equation:
\begin{eqnarray}
\label{wave} \frac{\partial^2 u}{\partial t^2}(t,x)&=& \Delta u
(t,x) +\dot
W(t,x), \quad t>0, x \in \bR^d \\
\nonumber
u(0, x)&=& 0, \quad x \in \bR^d \\
\nonumber \frac{\partial u}{\partial t}(0,x) &=& 0, \quad x \in
\bR^d.
\end{eqnarray}

Let $G_1$ be the fundamental solution of $u_{tt}-\Delta u=0$. It is
known that $G_1(t, \cdot)$ is a distribution in $\cS'(\bR^d)$ with
rapid decrease, and
\begin{equation}
\label{Fourier-G-wave} \cF
G_1(t,\cdot)(\xi)=\frac{\sin(t|\xi|)}{|\xi|},
\end{equation}
for any $\xi \in \bR^d,t>0,d \geq 1$ (see e.g. \cite{treves75}). In
particular,
\begin{eqnarray*}
G_1(t,x)&=&\frac{1}{2}1_{\{|x|<t\}}, \quad \mbox{if} \ d=1 \\
G_1(t,x)&=&\frac{1}{2 \pi}\frac{1}{\sqrt{t^2-|x|^2}}1_{\{|x|<t\}},
\quad \mbox{if} \ d=2 \\
G_1(t,x)&=&c_{d}\frac{1}{t}\sigma_t, \quad \mbox{if} \ d=3,
\end{eqnarray*}
where $\sigma_t$ denotes the surface measure on the 3-dimensional
sphere of radius $t$.

The solution of (\ref{wave}) is a square-integrable process
$u=\{u(t,x); t \geq 0, x \in \bR^d\}$ defined by:
$$u(t,x)=\int_{0}^{t} \int_{\bR^d}G_1(t-s,x-y)W(ds,dy).$$

By definition, $u(t,x)$ exists if and only if the stochastic
integral above is well-defined, i.e. $g_{tx}:=G_1(t-\cdot,x-\cdot)
\in \cH \cP$. In this case, $E|u(t,x)|^2 = \|g_{tx}\|_{\cH \cP}^2$.

The following theorem is the main result of this article.

\begin{theorem}
\label{wave-th} The solution $u=\{u(t,x); t \geq 0,x \in \bR^d\}$ of
(\ref{wave}) exists if and only if the measure $\mu$ satisfies
(\ref{wave-cond}). In this case, for all $p \geq 2$ and $T>0$
\begin{equation}\label{sup-L2-norm} \sup_{t \in [0,T]} \sup_{x \in
\bR^d} E|u(t,x)|^p<\infty,
\end{equation}
and the map $(t,x) \mapsto u(t,x)$ is continuous from $\bR_{+}
\times \bR^d$ into $L^2(\Omega)$.
\end{theorem}

\begin{example}
{\rm Let $f(x)=\gamma_{\alpha,d}|x|^{-(d-\alpha)}$ be the Riesz
kernel of order $\alpha \in (0,d)$. Then
$\mu(d\xi)=|\xi|^{-\alpha}d\xi$  and (\ref{wave-cond}) is equivalent
to $\alpha>d-2H-1$.}
\end{example}

\begin{example}
{\rm Let $f(x)=\gamma_{\alpha}\int_0^{\infty}w^{(\alpha-d)/2-1}
e^{-w}e^{-|x|^2/(4w)}dw$ be the Bessel kernel of order $\alpha>0$.
Then $\mu(d\xi)=(1+|\xi|^2)^{-\alpha/2}$ and (\ref{wave-cond}) is
equivalent to $\alpha>d-2H-1$. }
\end{example}

\begin{example}
{\rm Let $f(x)=\prod_{i=1}^{d}(\alpha_{H_i}|x_i|^{2H_i-2})$ be the covariance function of a fractional Brownian field with $H_i>1/2$ for all $i=1,\ldots,d$.
Then $\mu(d\xi)=\prod_{i=1}^{d}(c_{H_i}|\xi_i|^{-(2H_i-1)})$ and (\ref{wave-cond}) is
equivalent to $\sum_{i=1}^{d}(2H_i-1)>d-2H-1$. (This can be seen using the change of variables to the polar coordinates.)}
\end{example}

\begin{remark}
\label{rem-AB-equiv} {\rm Condition (\ref{wave-cond}) is equivalent
to
$$\int_{|\xi| \leq 1}\mu(d\xi)<\infty \quad \mbox{and} \quad \int_{|\xi| \geq
1}\frac{1}{|\xi|^{2H+1}}\mu(d\xi)<\infty.$$ }
\end{remark}

\noindent {\bf Proof of Theorem \ref{wave-th}:} Note that
$g_{tx}=G_1(t-\cdot,x-\cdot)$ satisfies conditions (i)-(iii) of
Theorem \ref{theorem-about-HP}. Hence, $g_{tx} \in \cH \cP$ (i.e.
the solution $u$ of (\ref{wave}) exists) if and only if $I_t<\infty$
for all $t>0$, where
$$I_t:=\alpha_H \int_{\bR^d} \int_0^t \int_0^t \cF
g_{tx}(u,\cdot)(\xi) \overline{\cF
g_{tx}(v,\cdot)(\xi)}|u-v|^{2H-2}dudv \mu(d\xi),$$ and
$E|u(t,x)|^2=\|g_{tx} \|_{\cH \cP}^2 =I_t$. Since $\cF
g_{tx}(u,\cdot)(\xi)=e^{-i\xi \cdot x} \overline{\cF
G_1(t-u,\cdot)(\xi)}$,
$$ I_t= \alpha_H
\int_{\bR^d} \int_0^t \int_0^t \cF G_1(u,\cdot)(\xi) \overline{\cF
G_1(v,\cdot)(\xi)}|u-v|^{2H-2}dudv \mu(d\xi).$$

\noindent Using (\ref{Fourier-G-wave}), we obtain:
$$I_t=\alpha_H \int_{\bR^d} \frac{\mu(d\xi)}{|\xi|^2}\int_0^t \int_0^t
\sin(u|\xi|)\sin(v|\xi|) |u-v|^{2H-2}dudv.$$

We split the integral $\mu(d\xi)$ into two parts, which correspond
to the regions $\{|\xi| \leq 1\}$ and $\{|\xi| \geq 1\}$. We denote
the respective integrals by $I_t^{(1)}$ and $I_t^{(2)}$. Since the
integrand is non-negative
$I_t<\infty$ if and only if $I_t^{(1)}<\infty$ and
$I_t^{(2)}<\infty$.

The fact that condition (\ref{wave-cond}) is sufficient for
$I_t<\infty$ follows by Proposition \ref{wave-prop-1} below. The
necessity follows by Proposition \ref{wave-prop-2} (using Remark
\ref{rem-AB-equiv}).

Relation (\ref{sup-L2-norm}) with $p=2$ follows from the estimates
obtained for $I_t=E|u(t,x)|^2$, using Proposition \ref{wave-prop-1}.
For arbitrary $p \geq 2$, we use the fact that $E|u(t,x)|^p \leq C_p
(E|u(t,x)|^2)^{p/2}$, since $u(t,x)$ is a Gaussian random variable.
The $L^2(\Omega)$-continuity is proved in Proposition
\ref{sol-cont-L2}. $\Box$

\vspace{3mm}

We begin with an auxiliary result. To simplify the notation, we
introduce the following functions: for $\lambda>0,\tau>0$, let
\begin{equation}
\label{def-ft-gt} f_t(\lambda,\tau)=\sin \tau \lambda t -\tau \sin
\lambda t, \quad g_t(\lambda,\tau)=\cos \tau \lambda t -\cos \lambda
t.
\end{equation}

\begin{lemma}
\label{lemma-sin-cos} For any $\lambda>0$ and $t>0$,
$$c_t^{(1)} \frac{\lambda^3}{1+\lambda^2} \leq \int_{\bR}\frac{1}{(\tau^2-1)^2}
[f_t^2(\lambda,\tau) +g_t^2(\lambda,\tau)]d\tau  \leq c_t^{(2)}
\frac{\lambda^3}{1+\lambda^2},$$ where $c_t^{(1)}=c_1 (t \wedge
t^3)$ and $c_t^{(2)}=c_2(t+t^3)$, for some positive constants
$c_1,c_2$.
\end{lemma}

\noindent {\bf Proof:} From the proof of Lemma \ref{H-norm-sin}, we
see that:
$$\frac{1}{(\tau^2-1)^2}[f_t^2(\lambda,\tau) +g_t^2(\lambda,\tau)]= |\cF_{0,\lambda t} \varphi(\tau)|^2,$$ where
$\varphi(x)=\sin x$. Using the Plancharel's identity
(\ref{Plancharel-lemma}), we obtain:
\begin{eqnarray*}
\lefteqn{\int_{\bR}\frac{1}{(\tau^2-1)^2}[f_t^2(\lambda,\tau)
+g_t^2(\lambda,\tau)]d\tau= \int_{\bR} |\cF_{0,\lambda t}
\varphi(\tau)|^2d\tau = }\\
& &  2\pi \int_0^{\lambda t}|\sin x|^2 dx=2\pi \lambda \int_0^t
|\sin \lambda s|^2 ds  = 2\pi \lambda^3 \int_0^t \frac{|\sin \lambda
s|^2}{\lambda^2} ds
\end{eqnarray*}

\noindent The result follows using (\ref{estimates}): (see e.g.
Lemma 6.1.2) of \cite{sanz-sole05})
$$c_t^{(1)}\frac{1}{1+\lambda^2} \leq \int_0^t
\frac{|\sin \lambda s|^2}{\lambda^2}ds \leq
c_t^{(2)}\frac{1}{1+\lambda^2}.$$ $\Box$

We denote by $N_t(\xi)$ the $\cH(0,t)$-norm of $u \mapsto \cF G_1(u,
\cdot)(\xi)$, i.e.
$$N_t(\xi)=\frac{\alpha_H}{|\xi|^2}\int_0^t \int_0^t \sin(u|\xi|)
\sin(v|\xi|) |u-v|^{2H-2}dudv.$$

\begin{proposition}
\label{wave-prop-1} For any $t>0, \xi \in \bR^d$
\begin{eqnarray*}
N_t(\xi) & \leq & C_{H}t^{2H+2} \left(\frac{1}{1+|\xi|^2}
\right)^{H+1/2}, \quad \mbox{if} \quad |\xi| \leq 1 \\
N_t(\xi) & \leq & c_{t,H}^{(3)} \left(\frac{1}{1+|\xi|^2}
\right)^{H+1/2}, \quad \mbox{if} \quad |\xi| \geq 1
\end{eqnarray*}
where $C_H= b_H^2 2^{H+1/2}/3$ and $c_{t,H}^{(3)}=c_H (
\frac{C}{1-H} + c_t^{(2)}) 2^{3H-1/2}$. Here $c_t^{(2)}$ is the
constant given by Lemma \ref{lemma-sin-cos}.
\end{proposition}

\noindent {\bf Proof:} a) Suppose that $|\xi| \leq 1$. We use the
fact that $\|\varphi \|_{\cH(0,t)}^2 \leq b_H^2 \|\varphi
\|_{L^{1/H}(0,t)}^2 \leq b_H^2 t^{2H-1} \|\varphi \|_{L^2(0,t)}^2$
for any $\varphi \in L^2(0,t)$, and $|\sin x| \leq x$ for any $x>0$.
Hence,
\begin{eqnarray*}
N_t(\xi) & \leq & b_H^2 t^{2H-1} \frac{1}{|\xi|^2}\int_0^t \sin^2
(u|\xi|)du \leq
 b_H^2 t^{2H-1}  \int_0^t u^2 du \\
 &=&  b_H^2 t^{2H-1} \frac{t^3}{3} \leq \frac{1}{3} b_H^2 t^{2H+2}
 2^{H+1/2}
 \left(\frac{1}{1+|\xi|^2} \right)^{H+1/2},
\end{eqnarray*}
where for the last inequality we used the fact that $\frac{1}{2}
\leq \frac{1}{1+|\xi|^2}$ if $|\xi| \leq 1$.

b) Suppose that $|\xi| \geq 1$. Using the change of variable
$u'=u|\xi|$, $v'=v|\xi|$,
\begin{eqnarray*}
N_t(\xi) &=&  \frac{\alpha_H }{|\xi|^{2H+2}}
\int_0^{t|\xi|}\int_0^{t|\xi|} \sin(u')\sin(v')|u'-v'|^{2H-2}dudv
 \\
&=&  \frac{1}{|\xi|^{2H+2}} \|\sin(\cdot)\|_{\cH(0,t|\xi|)}^2.
\end{eqnarray*}

Using the expression of the $\cH(0,t|\xi|)$-norm of $\sin(\cdot)$
given by Lemma \ref{H-norm-sin},
 we obtain:
\begin{equation}
\label{expression-N} N_t(\xi)=\frac{c_{H} }{|\xi|^{2H+2}} \int_{\bR}
\frac{|\tau|^{-(2H-1)}}{(\tau^2-1)^2} [f_t^2(|\xi|,\tau)+
g_t^2(|\xi|,\tau)]d\tau.
\end{equation}
We split the integral into the regions $|\tau| \leq 1/2$ and $|\tau|
\geq 1/2$, and we denote the two integrals by $N_{t}^{(1)}(\xi)$ and
$N_{t}^{(2)}(\xi)$.

Since $|f_t(\lambda,\tau)|\leq 1+|\tau|$ and $|g_t(\lambda,\tau)|
\leq 2$ for any $\lambda>0,\tau>0$, we have:
\begin{eqnarray*}
N_{t}^{(1)}(\xi) & \leq & c_H \frac{1}{|\xi|^{2H+2}} \int_{|\tau|
\leq 1/2} \frac{|\tau|^{-(2H-1)}}{(1-\tau^2)^2} [(1+|\tau|)^2+4]
d\tau \\
& \leq & c_H \frac{1}{|\xi|^{2H+1}}
\int_{|\tau| \leq 1/2} C|\tau|^{-(2H-1)} d\tau\\
&=& C
\frac{c_H}{1-H}\left(\frac{1}{2}\right)^{2-2H}\frac{1}{|\xi|^{2H+1}}.
\end{eqnarray*}
We used the fact that $|\xi|^{2H+2} \geq |\xi|^{2H+1}$ if $|\xi|
\geq 1$, and $\frac{1}{(1-\tau^2)^2}[(1+|\tau|)^2+4] \leq
\frac{1}{(3/4)^2}[(3/2)^2+4]=C$ if $|\tau| \leq 1/2$.

Using the fact that $|\tau|^{-(2H-1)} \leq (\frac{1}{2})^{-(2H-1)}$
if $|\tau| \geq \frac{1}{2}$, Lemma \ref{lemma-sin-cos}, and the
fact that $|\xi|^2/(1+|\xi|^2) \leq 1$, we obtain:
\begin{eqnarray*}
N_{t}^{(2)}(\xi) & \leq & \frac{c_H}{2^{-(2H-1)}}
\frac{1}{|\xi|^{2H+2}}\int_{|\tau| \geq 1/2}
\frac{1}{(\tau^2-1)^2}[f_{t}^2(|\xi|,\tau)+g_{t}^2(|\xi|,\tau)]d\tau \\
& \leq & \frac{c_H}{2^{-(2H-1)}}  \frac{1}{|\xi|^{2H+2}}\int_{\bR}
\frac{1}{(\tau^2-1)^2}[f_{t}^2(|\xi|,\tau)+g_{t}^2(|\xi|,\tau)]d\tau\\
& \leq & \frac{c_H}{2^{-(2H-1)}} c_{t}^{(2)}
\frac{1}{|\xi|^{2H+2}} \cdot |\xi| \frac{|\xi|^2}{1+|\xi|^2} \\
& \leq & \frac{c_H}{2^{-(2H-1)}} c_{t}^{(2)} \frac{1}{|\xi|^{2H+1}}.
\end{eqnarray*}
$\Box$

\begin{proposition}
\label{wave-prop-2} a) If $I_t^{(1)}<\infty$ for $t=1$, then
$\int_{|\xi| \leq 1}\mu(d\xi)<\infty$.

b) Let $l \geq 1$ be the integer from (\ref{mu-tempered}) and
$m=2l-2$. For any $t>0$,
\begin{equation}
\label{lower-estimate-It2} \int_{|\xi| \geq 1}
\frac{\mu(d\xi)}{|\xi|^{2H+1}} \leq a_{H,t}
(\sum_{i=0}^{m}b_t^{i})I_t^{(2)}+b_t^{m+1} \int_{|\xi| \geq
1}\frac{\mu(d\xi)}{|\xi|^{2H+2+m}},
\end{equation}
 where
$a_{H,t}=2^{2H} /(c_H c_{t}^{(1)})$, $b_t=2C /c_{t}^{(1)}$ and
$c_t^{(1)}$ is the constant of Lemma \ref{lemma-sin-cos}.

In particular, if $I_t^{(2)}<\infty$ for some $t>0$, then
$\int_{|\xi| \geq 1} |\xi|^{-(2H+1)} \mu(d\xi)<\infty$.
\end{proposition}

\noindent {\bf Proof:} a) Using the fact that $\sin x/x \geq \sin1$
for all $x \in [0,1]$, we have:
\begin{eqnarray*}
I_1^{(1)}&=&\int_{|\xi| \leq 1}\frac{\mu(d\xi)}{|\xi|^2} \int_0^1
\int_0^1 \sin (u|\xi|)\sin(v|\xi|)|u-v|^{2H-2}du dv \\
& \geq & \sin^2 1 \int_{|\xi| \leq 1}\mu(d\xi)\int_0^1 \int_0^1
uv|u-v|^{2H-2}du dv. 
\end{eqnarray*}

b) According to (\ref{expression-N}),
\begin{equation}
\label{new-expr-It2} I_t^{(2)}=c_{H} \int_{|\xi| \geq 1}
\frac{\mu(d\xi)}{|\xi|^{2H+2}} \int_{\bR}
\frac{|\tau|^{-(2H-1)}}{(\tau^2-1)^2} [f_t^2(|\xi|,\tau)+
g_t^2(|\xi|,\tau)]d\tau.
\end{equation}

 For any $k \in \{-1,0, \ldots,m\}$, let
$$I(k):=\int_{|\xi| \geq 1}\frac{1}{|\xi|^{2H+2+k}} \mu(d\xi).$$ By
(\ref{mu-tempered}), $I(m)= \int_{|\xi| \geq 1}|\xi|^{-(2H+2+m)}
\mu(d\xi)\leq \int_{|\xi| \geq 1}|\xi|^{-2l}\mu(d\xi)<\infty$.

We will prove that the integrals $I(k)$ satisfy a certain recursive
relation. By reverse induction, this will imply that all integrals
$I(k)$ with $k \in \{-1,0,\ldots,m\}$ are finite. For this, for $k
\in \{0,1 \ldots,m\}$, we let
\begin{equation}
 \label{claim2}
A_t(k):=\int_{|\xi| \geq
1}\frac{\mu(d\xi)}{|\xi|^{2H+2+k}}\int_{\bR}
\frac{1}{(\tau^2-1)^2}[f_t^2(|\xi|,\tau)+g_t^2(|\xi|,\tau)]d\tau.
\end{equation}

We consider separately the regions $\{|\tau| \leq 2\}$ and $\{|\tau|
\geq 2\}$. For the region $\{|\tau| \leq 2\}$, we use the expression
(\ref{new-expr-It2}) of $I_t^{(2)}$. Using the fact that
$|\xi|^{2H+2+k} \geq |\xi|^{2H+2}$ ({\em since $k \geq 0$}), and
$|\tau|^{-(2H-1)} \geq 2^{-(2H-1)}$ if $|\tau| \leq 2$, we obtain:
\begin{eqnarray*}
A_t'(k)& := &\int_{|\xi| \geq 1}\frac{\mu(d\xi)}{|\xi|^{2H+2+k}}
\int_{|\tau| \leq 2}\frac{1}{(\tau^2-1)^2}
[f_t^2(|\xi|,\tau)+g_t^2(|\xi|,\tau)]d\tau  \\
& \leq & 2^{2H-1} \int_{|\xi| \geq 1}\frac{\mu(d\xi)}{|\xi|^{2H+2}}
\int_{|\tau| \leq 2}\frac{|\tau|^{-(2H-1)}}{(\tau^2-1)^2}
[f_t^2(|\xi|,\tau)+g_t^2(|\xi|,\tau)]d\tau \\
& \leq &  2^{2H-1}\frac{1}{c_{H}}I_{t}^{(2)}, \quad \mbox{by}
(\ref{new-expr-It2}).
\end{eqnarray*}

For the region $\{|\tau| \geq 2\}$, we use the fact
$|f_t(\lambda,\tau)|\leq 1+|\tau|$ and $|g_t(\lambda,\tau)|\leq 2$
for all $\lambda>0,\tau>0$. Hence,
\begin{eqnarray*}
A_t''(k) &:=& \int_{|\xi| \geq 1}\frac{\mu(d\xi)}{|\xi|^{2H+2+k}}
\int_{|\tau| \geq 2}
\frac{1}{(\tau^2-1)^2}[f_t^2(|\xi|,\tau)+g_t^2(|\xi|,\tau)]d\tau
 \\
 & \leq & \int_{|\xi| \geq 1}\frac{\mu(d\xi)}{|\xi|^{2H+2+k}}
\int_{|\tau| \geq 2} \frac{1}{(\tau^2-1)^2}[(1+|\tau|)^2+4]d\tau
=CI(k).
\end{eqnarray*}
Hence, for any $k \in \{0,1,\ldots,m\}$
$$A_t(k) \leq 2^{2H-1}\frac{1}{c_{H}}I_{t}^{(2)}+CI(k).$$

Using Lemma \ref{lemma-sin-cos}, and the fact that
$\frac{|\xi|^2}{1+|\xi|^2} \geq \frac{1}{2}$ if $|\xi| \geq 1$, we
obtain: $$A_t(k) \geq c_t^{(1)} \int_{|\xi| \geq
1}\frac{\mu(d\xi)}{|\xi|^{2H+2+k}} \cdot \frac{|\xi|^3}{1+|\xi|^2}
\geq \frac{1}{2} c_t^{(1)} I(k-1),$$ for all $k \in
\{0,1,\ldots,m\}$. From the last two relations, we conclude that:
\begin{equation}
\label{recursion} \frac{1}{2} c_t^{(1)} I(k-1) \leq
 2^{2H-1}\frac{1}{c_{H}}I_{t}^{(2)}+CI(k), \quad \forall k \in
 \{0,1,\ldots,m\},
 \end{equation}
or equivalently, $I(k-1) \leq a_{H,t}I_t^{(2)}+b_t I(k)$, for all $k
\in \{0,1,\ldots,m\}$. Relation (\ref{lower-estimate-It2}) follows
by recursion. $\Box$

\begin{remark}
{\rm In the previous argument, the recursion relation
(\ref{recursion}) uses the fact that $k$ is {\em non-negative} (see
the estimate of $A_t'(k)$). Therefore, the ``last'' index $k$ for
which this relation remains true (counting downwards from $m$) is
$k=0$, leading us to the conclusion that $\int_{|\xi| \geq
1}|\xi|^{-(2H+1)} \mu(d\xi)<\infty$, if $I_t^{(2)}<\infty$. }
\end{remark}

The next result shows that the map $(t,x) \to u(t,x)$  from $\bR_{+}
\times \bR^d$ into $L^2(\Omega)$ is continuous.

\begin{proposition}
\label{sol-cont-L2} Suppose that (\ref{wave-cond}) holds, and let
$u=\{u(t,x), t \geq 0, x \in \bR^d\}$ be the solution of
(\ref{wave}). For any $t \geq 0$,
\begin{equation} \label{cont-t} E|u(t+h,x)-u(t,x)|^2 \to 0 \quad
\mbox{as} \ |h| \to 0, \quad \mbox{uniformly in} \ x \in
\bR^d\end{equation} and
\begin{equation}\label{cont-x}E|u(t,x)-u(t,y)|^2 \to 0 \quad
\mbox{as} \quad |x-y| \to 0.
\end{equation}
\end{proposition}

\noindent {\bf Proof:} We use the same argument as in Lemma 19 of
\cite{dalang99} (see also the erratum to \cite{dalang99}). We first
show (\ref{cont-t}).

Suppose that $h>0$. Splitting the interval $[0,t+h]$ into the
intervals $[0,t]$ and $[t,t+h]$, and using the inequality $|a+b|^2
\leq 2(a^2+b^2)$, we obtain:
\begin{eqnarray*}
E|u(t+h,x)-u(t,x)|^2 & \leq & 2 \{\|(g_{t+h,x}-g_{tx})1_{[0,t]} \|_{\cH
\cP}^2 +
 \|g_{t+h,x}1_{[t,t+h]} \|_{\cH \cP}^2\}\\
 &=:& 2[E_{1,t}(h)+E_2(h)].
 \end{eqnarray*}

Since  $\cF (g_{t+h,x}-g_{tx})(u,\cdot)(\xi)=e^{-i \xi \cdot x}
\overline{\cF G_1(t+h-u,\cdot)(\xi)-\cF G_1(t-u,\cdot)(\xi)}$,
\begin{eqnarray*}
E_{1,t}(h) &= & \alpha_H \int_{\bR^d} \mu(d\xi )\int_0^t \int_0^t dv dv
|u-v|^{2H-2}\cF (g_{t+h,x}-g_{tx})(u,\cdot)(\xi) \\
& & \overline{\cF
(g_{t+h,x}-g_{tx})(v,\cdot)(\xi)}\\
 &=& \alpha_H \int_{\bR^d} \mu(d\xi)\int_0^t \int_0^t du dv |u-v|^{2H-2}
[\cF G_1(u+h,\cdot)(\xi)-\cF G_1(u,\cdot)(\xi)] \\
& & \overline{\cF G_1(v+h,\cdot)(\xi)-\cF G_1(v,\cdot)(\xi)} \\
&=& \int_{\bR^d}\frac{\mu(d\xi)}{|\xi|^2} k_t(h,|\xi|),
\end{eqnarray*}
 where
\begin{eqnarray*}
k_t(h,|\xi|)&=&\alpha_H \int_0^t \int_0^t (\sin ((u+h)|\xi|) -\sin
(u|\xi|)) (\sin ((v+h)|\xi|) -\sin (v|\xi|)) \\
& & |u-v|^{2H-2}du dv =\| \sin ((\cdot \ +h)|\xi|) -\sin (\cdot \
|\xi|)\|_{\cH(0,t)}^2.
\end{eqnarray*}

\noindent By the Bounded Convergence Theorem, $\lim_{h \downarrow
0}k_t(h,|\xi|)=0$, for any $\xi \in \bR^d$.

The fact that $E_{1,t}(h) \to 0$ as $h \downarrow 0$ will follow from
the Dominated Convergence Theorem, once we prove that:
\begin{equation}
\label{bound-of-k} k_t(h,|\xi|) \leq k_t(|\xi|), \quad \forall h \in
[0,1], \forall \xi \in \bR^d, \quad \mbox{and} \quad
\int_{\bR^d}\frac{\mu(d\xi)}{|\xi|^2} k_t(|\xi|)<\infty.
\end{equation}

When $|\xi| \leq 1$, using the same argument as in Proposition
\ref{wave-prop-1}, we get:
\begin{eqnarray*}
k_t(h,|\xi|) & \leq & b_H^2 t^{2H-1} \| \sin ((\cdot \ +h)|\xi|)
-\sin
(\cdot \ |\xi|)\|_{L^2(0,t)}^2 \\
& \leq & 2b_H^2 t^{2H-1} \left(\int_0^t \sin^2
((u+h)|\xi|)du+\int_0^t \sin^2 (u|\xi|)du \right) \\
& \leq & 2 b_H^2 t^{2H-1}|\xi|^2\left(\int_0^t 2(u^2+1)du+\int_0^t
u^2 du \right)=:k_t(|\xi|).
\end{eqnarray*}

Suppose that $|\xi| \geq 1$. We use the fact that:
$$k_t(h,|\xi|)  \leq 2 (\| \sin ((\cdot \ +h)|\xi|)\|_{\cH(0,t)}^2+
\|\sin (\cdot \ |\xi|)\|_{\cH(0,t)}^2).$$

Using the change of variables $u'=(u+h)|\xi|, v'=(v+h)|\xi|$, and
(\ref{lemmaA1}) (Appendix A) we obtain:
\begin{eqnarray*}
\| \sin ((\cdot \ +h)|\xi|)\|_{\cH(0,t)}^2 &=& \alpha_H \int_0^t
\int_0^t  \sin ((u+h)|\xi|) \ \sin ((v+h)|\xi|) \ |u-v|^{2H-2}du dv \\
&=& \frac{\alpha_H}{|\xi|^{2H}}  \int_{h|\xi|}^{(t+h)|\xi|}
\int_{h|\xi|}^{(t+h)|\xi|} \sin (u') \sin (v') |u'-v'|^{2H-2}du' dv'
 \\
&=& \frac{c_H}{|\xi|^{2H}} \int_{\bR}  |\cF_{h|\xi|, (t+h)|\xi|}
\varphi (\tau)|^2 |\tau|^{-(2H-1)} d\tau,
\end{eqnarray*}
where $\varphi(t)=\sin t$. Note that the square of the real part of
$\cF_{h|\xi|, (t+h)|\xi|} \varphi(\tau)$ is:
$$\left|\int_{h|\xi|}^{(t+h)|\xi|} \cos \tau t \sin t dt\right|^2 \leq
2 \left|\int_{0}^{(t+h)|\xi|} \cos \tau t \sin t dt\right|^2 +2
\left|\int_{0}^{h|\xi|} \cos \tau t \sin t dt\right|^2,$$ and the
square of the imaginary part of $\cF_{h|\xi|, (t+h)|\xi|}
\varphi(\tau)$ is:
$$\left|\int_{h|\xi|}^{(t+h)|\xi|} \sin \tau t \sin t dt\right|^2 \leq
2 \left|\int_{0}^{(t+h)|\xi|} \sin \tau t \sin t dt\right|^2 +2
\left|\int_{0}^{h|\xi|} \sin \tau t \sin t dt\right|^2.$$

We now use the following fact (see Appendix B): for any $T>0$
$$\left| \int_0^T \cos \tau t \sin t dt \right|^2+
\left| \int_0^T \sin \tau t \sin t dt
\right|^2=\frac{1}{(\tau^2-1)^2}[(\sin \tau T-\tau \sin T)^2+(\cos
\tau T-\cos T)^2].$$

From here, it follows that $k_t(h,|\xi|)$ is bounded by:
$$\frac{2c_H}{|\xi|^{2H}}
\int_{\bR}\frac{|\tau|^{-(2H-1)}}{(\tau^2-1)^2}[f_{t+h}^2(|\xi|,\tau)+
g_{t+h}^2(|\xi|,\tau)+f_{h}^2(|\xi|,\tau)+g_{h}^2(|\xi|,\tau)+
f_{t}^2(|\xi|,\tau)+g_{t}^2(|\xi|,\tau)]d\tau,$$ where
$f_t(\lambda,\tau)$ and $g_{t}(\lambda,\tau)$ are defined by
(\ref{def-ft-gt}). The argument of Proposition \ref{wave-prop-1}
shows that for any $t>0$ and $|\xi| \geq 1$
$$\int_{\bR}
\frac{|\tau|^{-(2H-1)}
}{(\tau^2-1)^2}[f_t^2(|\xi|,\tau)+g_t^2(|\xi|,\tau)] d\tau \leq
c_{t,H}^{(4)}\frac{|\xi|^3}{1+|\xi|^2} \leq c_{t,H}^{(4)}|\xi|,$$
where $c_{t,H}^{(4)}=\frac{2C}{1-H}\left(\frac{1}{2} \right)^{2-2H}
+ \left(\frac{1}{2}\right)^{-(2H-1)}c_t^{(2)}$. Since
$c_{t,H}^{(4)}$ is non-decreasing in $t$ and $h \in [0,1]$,
$k_t(h,|\xi|)$ is bounded by
$$\frac{2c_H}{|\xi|^{2H-1}} (c_{t+h,H}^{(4)}+c_{h,H}^{(4)}+c_{t,H}^{(4)})
\leq \frac{2c_H}{|\xi|^{2H-1}}
(c_{t+1,H}^{(4)}+c_{1,H}^{(4)}+c_{t,H}^{(4)}):=k_t(|\xi|).$$ This
concludes the proof of (\ref{bound-of-k}).

A similar argument shows that $E_2(h) \to 0$ as $h \downarrow 0$,
since
\begin{eqnarray*}
E_2(h) &=& \alpha_H \int_{\bR^d} \int_{t}^{t+h} \int_t^{t+h} \cF G_1
(t+h-u, \cdot)(\xi) \overline{\cF G_1 (t+h-v,
\cdot)(\xi)}|u-v|^{2H-2}du dv \mu(d\xi) \\
&=& \alpha_H \int_{\bR^d}\frac{\mu(d\xi)}{|\xi|^2} \int_0^h \int_0^h
\sin (u|\xi|) \sin (v|\xi|) |u-v|^{2H-2}du dv.
\end{eqnarray*}

The case $h<0$ is treated similarly. Using the same argument as
above, it follows that for any $h>0$, $E|u(t-h,x)-u(t,x)|^2 \leq 2
(E_{1,t}'(h)+E_2(h))$, where
$$E_{1,t}'(h)=\int_{\bR^{d}} \frac{\mu(d\xi)}{|\xi|^2} k_t'(h,|\xi|), \
\mbox{and} \ k_t'(h,|\xi|)=\| \sin (\cdot \ |\xi|)-\sin ( (\cdot \ -h) |\xi|)  \|_{\cH(h,t)}^2.$$


 To prove (\ref{cont-x}), note that
\begin{eqnarray*}
\lefteqn{E|u(t,x)-u(t,y)|^2=\|g_{tx}-g_{ty}\|_{\cH \cP}^2 = }\\
& & \alpha_H \int_{\bR^d} \int_{0}^{t} \int_{0}^{t} \cF
(g_{tx}-g_{ty}) (u, \cdot)(\xi) \overline{\cF (g_{tx}-g_{ty})(v,
\cdot)(\xi)}|u-v|^{2H-2}du dv \mu(d\xi)= \\
& & \alpha_H \int_{\bR^d} |e^{-i \xi \cdot x}-e^{-i \xi \cdot
y}|^2\int_{0}^{t} \int_{0}^{t} \cF G_1(u, \cdot)(\xi) \overline{\cF
G_1(v, \cdot)(\xi)}|u-v|^{2H-2}du dv \mu(d\xi),
\end{eqnarray*}
which converges to $0$ as $|x-y| \to 0$, by the Dominated
Convergence Theorem. $\Box$

\begin{example}
{\rm There exists an interesting connection between the solution of
the wave equation with fractional noise in time and Riesz covariance
in space and the odd and even parts of the fBm. Indeed, if $f$ be
the Riesz kernel of order $\alpha \in (0,d)$, then
\begin{eqnarray*}
I_t&=& \alpha_H \int_{\bR^d}d\xi |\xi|^{-\alpha-2H-2}\int_{\bR}
\frac{|\tau|^{-(2H-1)}}{(\tau^2-1)^2}[f_t^2(|\xi|,\tau)+g_t^2(|\xi|,\tau)]
d\tau \\
&=& 2\alpha_H c_d \int_{\bR} \frac{|\tau|^{-(2H-1)}}{(\tau^2-1)^2}
\left(\int_0^{\infty} \frac{(\sin \tau \lambda t -\tau \sin \lambda
t)^2}{\lambda^2}\lambda^{-\theta}d\lambda+\int_0^{\infty}
\frac{(\cos \tau \lambda t -\cos \lambda
t)^2}{\lambda^2}\lambda^{-\theta}d\lambda \right),
\end{eqnarray*}
where $\theta=\alpha+1-d+2H >0$ under (\ref{wave-cond}). If
$\theta<1$, the two integrals $d\lambda$ can be expressed in terms
of the covariance functions of the odd and  even parts of the fBm
(see \cite{DZ}).}
\end{example}

\section{The heat equation}

In this section, we consider the the heat equation with additive
noise:
\begin{eqnarray}
\label{heat} \frac{\partial u}{\partial t}(t,x)&=&\frac{1}{2} \Delta
u(t,x)+\dot
W(t,x), \quad t>0, x \in \bR^d \\
\nonumber u(0,x)&=& 0, \quad x \in \bR^d.
\end{eqnarray}

Equation (\ref{heat}) was treated in \cite{balan-tudor08}, in the
case of particular covariance kernels $f$. We give here an unitary
approach which covers the case of any covariance kernel $f$, which
satisfies (\ref{heat-cond}).

The case of the heat equation is actually much simpler than the case of the
wave equation, since both the fundamental solution $G$ and its
Fourier transform are non-negative functions.

More precisely, let $G_2$ be the fundamental solution of
$u_t-\frac{1}{2}\Delta u=0$. Then
$$G_2(t,x)=\frac{1}{(2\pi t)^{d/2}} \exp \left(-\frac{|x|^2}{2t}\right), \quad t >0,x \in
\bR^d$$ and
\begin{equation}
\label{Fourier-heat} \cF G_2(t,\cdot)(\xi)=\exp
\left(-\frac{t|\xi|^2}{2} \right), \quad t>0,\xi \in \bR^d.
\end{equation}

We will prove the following result.

\begin{theorem}
\label{heat-th} The solution $u=\{u(t,x),t \geq 0, x \in \bR^d\}$ of
(\ref{heat}) exists if and only if the measure $\mu$ satisfies
(\ref{heat-cond}). In this case, (\ref{sup-L2-norm}) holds for all
$p \geq 2$ and $T>0$, and the solution is $L^2(\Omega)$-continuous.
\end{theorem}

\begin{remark}
{\rm (i) When $f$ is the Riesz kernel of order $\alpha$, or the
Bessel kernel of order $\alpha$, condition (\ref{heat-cond}) is
equivalent to $\alpha>d-4H$. When $f$ is the covariance function of the fractional Brownian field with
$H_i>1/2$ for all $i=1,\ldots,d$, condition (\ref{heat-cond}) is
equivalent to $\sum_{i=1}^{d}(2H_i-1)>d-4H$. Note
that this condition is weaker than the condition given in \cite{oksendal-zhang01}.

(ii) In Theorem 2.1 of the Erratum to \cite{balan-tudor08} it has
been proven that condition (\ref{heat-cond}) implies that
$\|g_{tx}\|_{\cH \cP}<\infty$ for any $t\geq 0$ and $x\in \bR ^{d}$.
}
\end{remark}


\noindent {\bf Proof of Theorem \ref{heat-th}:} Note that
$g_{tx}=G_2(t-\cdot,x-\cdot)$ is non-negative. Hence, $g_{tx} \in
\cH \cP$ 
if and only if $g_{tx} \in |\cH \cP|$. This is equivalent to saying
that $J_t:=\|g_{tx}\|_{|\cH \cP|}^2<\infty$ for all $t>0$. Note that
\begin{eqnarray*}
J_t&=& \alpha_H  \int_0^t \int_0^t \int_{\bR^d} \int_{\bR^d}
 g_{tx}(u,y)
g_{tx}(v,z)f(y-z)|u-v|^{2H-2}dy dz dudv   \\
&=& \alpha_H \int_0^t \int_0^t  \int_{\bR^d} \cF
g_{tx}(u,\cdot)(\xi) \overline{\cF g_{tx}(v,\cdot)(\xi)}|u-v|^{2H-2}
\mu(d\xi) dudv
 \\
&=& \alpha_H   \int_0^t \int_0^t \int_{\bR^d} \cF
G_2(t-u,\cdot)(\xi) \overline{\cF G_2(t-v,\cdot)(\xi)}|u-v|^{2H-2}
 \mu(d\xi) dudv.
\end{eqnarray*}

Using (\ref{Fourier-heat}) and Fubini's theorem (whose application
is justified since the integrand is non-negative), we obtain:
$$J_t=\alpha_H \int_{\bR^d} \int_0^t \int_0^t \exp
\left(-\frac{u|\xi|^2}{2} \right) \exp \left(-\frac{v|\xi|^2}{2}
\right)  |u-v|^{2H-2}dudv \mu(d\xi).$$

The existence of the solution follows from Proposition
\ref{estimate-heat} below, which also gives estimates for
$J_t=E|u(t,x)|^2$ (and hence for $E|u(t,x)|^p$). The
$L^2(\Omega)$-continuity is given by Proposition \ref{cont-L2-heat}.
$\Box$

Let
$$A_t(\xi)=\alpha_H\int_0^t \int_0^t \exp
\left(-\frac{u|\xi|^2}{2} \right) \exp \left(-\frac{v|\xi|^2}{2}
\right)  |u-v|^{2H-2}dudv $$

The next result is similar to Lemma 6.1.1) of \cite{sanz-sole05}.

\begin{proposition}
\label{estimate-heat} For any $t>0, \xi \in \bR^d$,
$$\frac{1}{4} (t^{2H} \wedge 1)\left(\frac{1}{1+|\xi|^2}\right)^{2H}
\leq A_t(\xi) \leq C_H (t^{2H}+1)
\left(\frac{1}{1+|\xi|^2}\right)^{2H},$$ where $C_H=b_H^2
(4H)^{2H}$.
\end{proposition}

\noindent {\bf Proof:} Suppose that $|\xi| \leq 1$. Using the fact
that $\|\varphi \|_{\cH(0,t)}^2 \leq b_H^2 t^{2H-1}
\|\varphi\|_{L^2(0,t)}^2$ for all $\varphi \in L^2(0,t)$, $e^{-x}
\leq 1$ for any $x>0$, and $\frac{1}{2} \leq \frac{1}{1+|\xi|^2}$ if
$|\xi| \leq 1$,
$$A_t(\xi) \leq b_H^2 t^{2H-1} \int_0^t
\exp(-u|\xi|^2)du \leq b_H^2 t^{2H} \leq b_H^2 2^{2H} t^{2H}
\left(\frac{1}{1+|\xi|^2}\right)^{2H}.$$

Suppose that $|\xi| \geq 1$. Using the fact that $\|\varphi
\|_{\cH(0,t)}^2 \leq b_H^2 \|\varphi\|_{L^{1/H}(0,t)}^2$ for any
$\varphi \in L^{1/H}(0,t)$, $1-e^{-x} \leq 1$ for all $x>0$, and
$\frac{1}{|\xi|^2} \leq \frac{2}{1+|\xi|^2}$, we obtain:
\begin{eqnarray*}
A_t(\xi) & \leq & b_{H}^2 \left[\int_0^t
\exp\left(-\frac{u|\xi|^2}{2H} \right) du \right]^{2H} =  b_{H}^2
\left(\frac{2H}{|\xi|^2}\right)^{2H}
\left[1-\exp\left(-\frac{t|\xi|^2}{2H}\right) \right]^{2H} \\
& \leq &   b_{H}^2 (4H)^{2H}
\left(\frac{1}{1+|\xi|^{2}}\right)^{2H}.
\end{eqnarray*}
This proves the upper bound.

Next, we show the lower bound. Suppose that $t|\xi|^2 \leq 1$. For
any $u \in [0,t]$, $\frac{u|\xi|^2}{2} \leq \frac{t|\xi|^2}{2} \leq
\frac{1}{2}$. Using the fact that $e^{-x} \geq 1-x$ for all $x>0$,
we conclude that:
$$\exp\left(-\frac{u|\xi|^2}{2}\right) \geq 1-\frac{u|\xi|^2}{2}
\geq \frac{1}{2}, \quad \forall u \in [0,t].$$ Hence
$$A_t(\xi) \geq \alpha_H \left(\frac{1}{2} \right)^2
\int_0^t \int_0^t |u-v|^{2H-2}dudv=\frac{1}{4}t^{2H} \geq
\frac{1}{4}t^{2H} \left( \frac{1}{1+|\xi|^2}\right)^{2H}.$$ For the
last inequality, we used the fact that $1 \geq \frac{1}{1+|\xi|^2}$.

Suppose that $t|\xi|^2 \geq 1$. Using the change of variables
$u'=u|\xi|^2/2, v'=v|\xi|^2/2$, we obtain:
$$A_{t}(\xi)=\alpha_H \frac{2^{2H}}{|\xi|^{4H}}  \int_0^{t|\xi|^2/2}
\int_{0}^{t|\xi|^2/2} e^{-u'}e^{-v'}|u'-v'|^{2H-2}du' dv'.$$ Since
the integrand is
non-negative, 
\begin{eqnarray*}
A_t(\xi) & \geq & \alpha_H \frac{2^{2H}}{|\xi|^{4H}}  \int_0^{1/2}
\int_{0}^{1/2} e^{-u}e^{-v}|u-v|^{2H-2}dudv \\
& = & 2^{2H}\|e^{-u} \|_{\cH(0,1/2)}^2 \frac{1}{|\xi|^{4H}} \geq
2^{2H} \left(\frac{1}{2} \right)^{2H+2} \left(\frac{1}{1+|\xi|^2}
\right)^{2H},
\end{eqnarray*}
where for the last inequality we used the fact that
$\frac{1}{|\xi|^2} \geq \frac{1}{1+|\xi|^2}$, and $\|e^{-u}
\|_{\cH(0,1/2)}^2 \geq \left(\frac{1}{2} \right)^{2H+2}$. (This
follows since $e^{-u} \geq 1-u \geq \frac{1}{2}$ for all $u \in
[0,\frac{1}{2}]$.)

$\Box$

\begin{proposition}
\label{cont-L2-heat} Suppose that (\ref{heat-cond}) holds, and let
$u=\{u(t,x), t \geq 0, x \in \bR^d\}$ be the mild-sense solution of
(\ref{heat}). Then the map $(t,x) \to u(t,x)$  from $\bR_{+} \times
\bR^d$ into $L^2(\Omega)$ is continuous.
\end{proposition}

\noindent {\bf Proof:} The argument is similar to that of
Proposition \ref{sol-cont-L2}. In this case, if $h>0$,
$$E_{1,t}(h)=\int_{\bR^d}\mu(d\xi)k_t(h,|\xi|),$$
where
$$k_t(h,|\xi|)=\left\| \exp \left(-\frac{(\cdot +h)|\xi|^2}{2} \right)-\exp
\left(-\frac{ \cdot  \ |\xi|^2}{2} \right) \right\|_{\cH(0,t)}^2,$$
and
$$E_2(h)=\alpha_H\int_{\bR^d}\mu(d\xi)\int_0^h \int_0^h \exp
\left(-\frac{u|\xi|^2}{2} \right) \exp \left(-\frac{v|\xi|^2}{2}
\right)|u-v|^{2H-2}du dv .$$ We omit the details. $\Box$

\begin{remark}
{\rm We consider the operator $Lu=\partial_t
u-\sum_{i,j=1}^{d}a_{ij}\partial^2_{x_i x_j}u-\sum_{i=1}^{d}b_i
\partial_{x_i}u$. Let $G_3(t,x;s,y)$ be the fundamental
solution of $Lu=0$. We assume that:

(i) The functions $a_{ij},b_i:[0,T] \times \bR^d \to \bR$,
$i,j=1,\ldots,d$ are $\alpha/2$-H\"older continuous in $t$ and
$\alpha$-H\"older continuous in $x$, for some $\alpha \in (0,1)$.

(ii) There exist some $k,K>0$ such that for all $(t,x) \in [0,T]
\times \bR^d, \xi \in \bR^d$,
$$k|\xi|^2 \leq \sum_{i,j=1}^{d}a_{ij}(t,x)\xi_i \xi_j \leq
K|\xi|^2.$$

Under these assumptions, $G_3$ is a positive function defined on
$[0,T] \times \bR^d \times [0,T] \times \bR^d \cap \{(s,t); 0 \leq s
\leq t \leq T\}$, which satisfies: (see p. 376 of \cite{LSU68})
\begin{equation}
\label{ineq-G3} G_3(t,x;s,y) \leq c_1 (t-s)^{-d/2} \exp \left(-c_2
\frac{|x-y|^2}{t-s} \right):=G_2'(t-s,x-y).
\end{equation}

Since $G_2'(t,x)$ is essentially the same as the heat kernel
$G_2(t,x)$, the solution of $Lu(t,x)=\dot W(t,x)$ (with vanishing
initial conditions) exists, if the measure $\mu$ satisfies condition
(\ref{heat-cond}). }
\end{remark}

\appendix

\section{Some useful identities}

Recall that the Fourier transform of a function $\varphi \in
L^1(\bR)$ is defined by: $$\cF \varphi(\tau)=\int_{\bR}e^{-i\tau
x}\varphi(x)dx.$$

\noindent For an interval $(a,b) \subset \bR$, we define the
restricted Fourier transform of a function $\varphi \in L^1(a,b)$:
$$\cF_{a,b} \varphi(\tau):=\int_{a}^{b} e^{-i\tau x}\varphi(x)dx=
\cF (\varphi 1_{[a,b]})(\tau).$$

One can prove that $\cF \varphi \in L^2(\bR)$, for any $\varphi \in
L^1(\bR) \cap L^2(\bR).$ By the Plancharel's identity, for any
$\varphi,\psi \in L^1(\bR) \cap L^2(\bR)$, we have:
$$ \int_{\bR} \varphi(x)\psi(x)dx =(2\pi)^{-1}\int_{\bR}\cF \varphi(\tau)
\overline{\cF \psi(\tau)}d\xi.$$

In particular, for any $\varphi,\psi \in L^2(a,b)$, we have:
\begin{equation}
\label{Plancharel-lemma} \int_{a}^{b} \varphi(x)\psi(x)dx
=(2\pi)^{-1}\int_{\bR}\cF_{a,b} \varphi(\tau) \overline{\cF_{a,b}
\psi(\tau)}d\xi.
\end{equation}

\noindent (Consider $\tilde \varphi=\varphi 1_{[a,b]}$. Then $\tilde
\varphi \in L^1(\bR) \cap L^2(\bR)$ and $\cF \tilde
\varphi(\xi)=\cF_{a,b} \varphi(\xi)$.)

\vspace{2mm}

The proof of Theorem \ref{wave-th} uses, in an essential way, a
formula for the $\cH(0,T)$-norm of $\sin$ (developed in Appendix B),
which is in turn based on the following result. (This result can be
derived using for instance, the results of \cite{PT00}.)


\begin{lemma}
Let $H \in (\frac{1}{2},1)$. For any $\varphi,\psi \in L^1(\bR) \cap
L^2(\bR)$,
$$\alpha_H \int_{\bR} \int_{\bR} \varphi(u)
\psi(v)|u-v|^{2H-2}du dv=c_{H} \int_{\bR} \cF \varphi
(\tau)\overline{\cF \psi(\tau)}|\tau|^{-(2H-1)}d\tau,$$ where
$\alpha_H=H(2H-1)$ and $c_H=\Gamma(2H+1)\sin(\pi H)/(2\pi)$.

In particular, for any $\varphi,\psi \in L^2(a,b)$,
\begin{equation}
\label{lemmaA1} \alpha_H \int_a^b \int_a^b \varphi(u)
\psi(v)|u-v|^{2H-2}du dv=c_{H} \int_{\bR} \cF_{a,b}\varphi
(\tau)\overline{\cF_{a,b}\psi(\tau)}|\tau|^{-(2H-1)}d\tau.
\end{equation}
\end{lemma}

\section{The $\cH(0,T)$-norm of $\sin$}

\begin{lemma}
\label{H-norm-sin} Let $\varphi(t)=\sin t$, $t \in [0,T]$. Then
$$\|\varphi\|_{\cH(0,T)}^2= c_H \int_{\bR} \frac{(\sin \tau T-\tau \sin T)^2+(\cos\tau
T-\cos T)^2}{(\tau^2-1)^2} \ |\tau|^{-(2H-1)} d\tau,$$ where
$c_H=\Gamma(2H+1)\sin(\pi H)/(2\pi)$.
\end{lemma}

\noindent {\bf Proof}: By (\ref{lemmaA1}),
$$\|\varphi \|_{\cH(0,T)}^2=\alpha_H \int_0^T \int_0^T \varphi(u)
\varphi(v) |u-v|^{2H-2}du dv=c_H \int_{\bR}|\cF_{0,T} \varphi
(\tau)|^2 |\tau|^{-(2H-1)}d\tau.$$ Note that
$|\cF_{0,T}\varphi(\tau)|^2=\left|\int_{0}^{T}e^{-i\tau
t}\varphi(t)dt\right|^2=I_1^2+ J_1^2$, where
$$I_1={\rm Re}[\cF_{0,T} \varphi (\tau)]=\int_0^T \cos \tau t \sin t dt,
\quad J_1={\rm Im}[\cF_{0,T} \varphi (\tau)]=\int_0^T \sin \tau t
\sin t dt.$$

\noindent We calculate $I_1$ first. Using integration by parts, we
obtain: $$ I_1=1-\cos \tau T\cos T-\tau I_2,$$ where $I_2=\int_0^T
\sin \tau t \cos t dt$. On the other hand,
$$I_1+I_2=\int_0^T \sin [(\tau+1 )t]dt=\frac{1-\cos[(\tau+1)T]}{\tau+1}.$$

\noindent Solving for $I_1$ and $I_2$, we obtain:
$$I_1=\frac{1}{1-\tau^2}(1-\cos\tau T \cos T- \tau \sin\tau T \sin T).$$

\noindent Similarly, letting $J_2=\int_0^T \cos\tau t \cos t dt$,
obtain:
$$\tau J_2-J_1=\sin\tau T \cos T \quad \mbox{and} \quad
J_2-J_1=\frac{1}{\tau+1}\sin[(\tau+1)T].$$ Solving for $J_1$, we
obtain:
$$J_1=\frac{1}{1-\tau^2}(\tau \cos\tau T \sin T-\sin\tau T \cos T).$$

\noindent An elementary calculation shows that:
$$I_1^2+J_1^2 =\frac{1}{(1-\tau^2)^2} [(\sin \tau T-\tau \sin
T)^2+(\cos\tau T-\cos T)^2].$$ $\Box$

\begin{remark}
{\rm Let $B=(B_t)_{t \in \bR}$ a fBm of index $H$ (on the whole real
line). Let $B^o=(B_t^o)_{t \in \bR}$ and $B^e=(B_t^e)_{t \in \bR}$
be the odd and even parts of $B$ (see \cite{DZ}). $B^o$ and $B^e$
are independent centered Gaussian processes with
$B_t=B_t^{o}+B_{t}^{e}$, and
\begin{eqnarray*} E(B_t^o
B_s^o)&=&\langle 1_{(0,t)}, 1_{(0,s)} \rangle_{o}:= c_{H}\int_{\bR}
\frac{\sin \tau t \sin \tau s}{\tau^2} |\tau|^{-(2H-1)}d\tau \\
E(B_t^e B_s^e)&= &\langle 1_{(0,t)}, 1_{(0,s)} \rangle_{e}:=
c_{H}\int_{\bR} \frac{(1-\cos \tau t)(1- \cos \tau s)}{\tau^2}
|\tau|^{-(2H-1)}d\tau.
\end{eqnarray*}

In general, for $\varphi \in L^1(\bR) \cap L^2(\bR)$, we have:
$$E[B^o(\varphi)^2]=\|\varphi \|_{o}^2:=c_{H} \int_{\bR}
|{\rm Re}[\cF \varphi (\tau)]|^2 |\tau|^{-(2H-1)}d\tau$$
$$E[B^e(\varphi)^2]=\|\varphi \|_{e}^2:=c_{H} \int_{\bR}
|{\rm Im}[\cF \varphi (\tau)]|^2 |\tau|^{-(2H-1)}d\tau$$

In the proof of Lemma \ref{H-norm-sin}, $I_1=I_1(\tau)$ and
$J_1=J_1(\tau)$ are the real and imaginary parts of $\cF (\varphi
1_{[0,T]})(\tau)$, where $\varphi (t)=\sin t$. Hence
\begin{eqnarray*}
E \left(\int_0^T \sin t dB_t^o \right)^2=\|\sin (\cdot)
1_{[0,T]}\|_{o}^2 &=& c_H \int_{\bR} |I_1(\tau)|^2
|\tau|^{-(2H-1)}d\tau \\
E \left(\int_0^T \sin t dB_t^e \right)^2=\|\sin (\cdot)
1_{[0,T]}\|_{e}^2 &=& c_H \int_{\bR} |J_1(\tau)|^2
|\tau|^{-(2H-1)}d\tau.
\end{eqnarray*}

 }
\end{remark}

\vspace{3mm}

\footnotesize{{\em Acknowledgement.} The authors would like to thank
Professor Robert Dalang for the invitations to visit EPFL and the
useful discussions.

\normalsize{

\end{document}